\documentclass{article}

\usepackage{graphicx,amsmath,amssymb,color,float}

\newtheorem{Theorem}{\quad Theorem}[section]

\newtheorem{Lemma}[Theorem]{\quad Lemma}

\title{Bifurcations in Delayed Lotka-Volterra\\ Intraguild Predation Model}
\author{Juancho A. Collera\\ Department of Mathematics and Computer Science\\ University of the Philippines Baguio\\ \texttt{jacollera@up.edu.ph}}

\begin{document}
\maketitle

\begin{abstract}
Omnivory is defined as feeding on more than one trophic level. An example of this is the so-called intraguild predation (IG) which includes a predator and its prey that share a common resource. IG predation models are known to exhibit interesting dynamics including chaos. This work considers a three-species food web model with omnivory, where the interactions between the basal resource, the IG prey, and the IG predator are of Lotka-Volterra type. In the absence of predation, the basal resource follows a delayed logistic equation or popularly known as Hutchinson's equation. Conditions for the existence, stability, and bifurcations of all non-negative equilibrium solutions are given using the delay time as parameter. Results are illustrated using numerical bifurcation analysis. 
\end{abstract}

\vspace{0.2in}


\section{Introduction}
Omnivory, as defined in \cite{pimmlawton}, occurs when a population feeds on resources at more than one trophic level. For example, a species feeding on its prey's resource is called omnivorous. This particular type of tri-trophic community module is called intraguild (IG) predation \cite{polismyersholt}, and was shown to be quite common in nature \cite{arummarquet}. 
This three-species IG predation model includes top and intermediate predators termed as IG predator and IG prey, respectively, and a basal resource. The IG predator depends completely both on the IG prey and the basal resource for its sustenance, while the IG prey depends solely on the basal resource. 
A study of a model of IG predation with non-linear functional responses in \cite{mccannhastings} showed that omnivory stabilizes and enhances persistence of the three-species food web. 
However, a model of Lotka-Volterra type with linear functional responses considered in \cite{holtpolis} showed that IG predation could have a destabilizing effect, and a criterion for co-existence of all three species is that the IG prey must be superior than the IG predator in competing for the shared basal resource. Lotka-Volterra IG predation models are known to exhibit interesting dynamics such as limit cycles \cite{holt}, bistabity \cite{hsuruanyang}, and chaos \cite{namba2005,namba2008}.\\

In this paper, we consider a Lotka-Volterra IG predation model where, in the absence of the IG predator and the IG prey, the basal resource grows according to the delayed logistic equation or more commonly known as the Hutchinson's equation \cite{arinonato}. It should be noted that there are several other ways of formulating the delayed logistic equation such as in \cite{arinoalternative}, and in this paper, we use the so-called \emph{classical} delayed logistic equation. Our delayed model generalizes the Lotka-Volterra IG predation models of \cite{holtpolis} and \cite{namba2005} which are basically systems of ODEs.\\

This paper is organized as follows. In Section 2, we  introduce the delayed Lotka-Volterra IG predation model and discuss the existence of its equilibrium solutions. In Section 3, we give the main results of this paper. These are the conditions for stability and bifurcations of all non-negative equilibria using the delay time as parameter. In Section 4, we use numerical continuation and bifurcation analysis to illustrate our results on the effects of the delay time to our Lotka-Volterra IG predation model. We then end by giving a summary and the thoughts of this paper.

\section{The Model}
We consider the following delayed Lotka-Volterra IG predation model
\begin{eqnarray}
\dot{x}(t)&=&\left[a_{0}-a_{1}x(t-\tau)-a_{2}y(t)-a_{3}z(t)\right]x(t),\nonumber\\ 
\dot{y}(t)&=&\left[-b_{0}+b_{1}x(t)-b_{3}z(t)\right]y(t),\label{eq:model}\\
\dot{z}(t)&=&\left[-c_{0}+c_{1}x(t)+c_{2}y(t)\right]z(t),\nonumber
\end{eqnarray}
where $x(t)$, $y(t)$, and $z(t)$ are the densities at time $t$ of the basal resource, IG prey, and IG predator, respectively. In the absence of the IG predator and the IG prey, the basal resource grows according to the delay logistic equation \cite{arinonato, halsmith}. When $\tau=0$, we recover the ODEs model in \cite{namba2005}. We refer to \cite{namba2005} for the description of the rest of the parameters in (\ref{eq:model}). Moreover, we use an initial condition $(x(t),y(t),z(t))=(x_0,y_0,z_0)$ for $t\in[-\tau,0]$ and where $x_0,y_0,z_0\ge0$.\\

An important characteristic of the IG predation is that is it a mixture of community modules such as competition and predation \cite{holt, holtpolis, polismyersholt}. When $a_3=c_1=0$, we obtain a food chain while if $b_3=c_2=0$, we obtain exploitative competition (or shared resources) where two predators, in our case the IG predator and the IG prey, share a common resource and the IG predator does not feed on the IG prey. Meanwhile, if $a_2=b_1=0$, we obtain apparent competition (or shared predation), which in our model, the IG predator feeds on both the  IG prey and the basal resource but there is no predation on the basal resource by the IG prey. In this case, the IG prey will go extinct.
\subsection{Equilibrium Solutions}
System (\ref{eq:model}) has five possible non-negative equilibrium solutions: 
$E_0=(0,0,0)$, 
$E_1=\left(K,0,0\right)$ where $K=a_0/a_1$, 
$E_2=\left(A,B,0\right)$ where $A=b_0/b_1$ and $B=(a_{0}b_{1}-a_{1}b_{0})/a_{2}b_{1}$, 
$E_3=\left(C,0,D\right)$ where $C=c_0/c_1$ and $D=(a_{0}c_{1}-a_{1}c_{0})/a_{3}c_{1}$, 
and the positive equilibrium solution 
$E_4=(P/S,Q/S,R/S)$ 
where 
\begin{eqnarray*}
P&=&a_{0}b_{3}c_{2}-a_{2}b_{3}c_{0}+a_{3}b_{0}c_{2},\\
Q&=&-a_{0}b_{3}c_{1}+a_{1}b_{3}c_{0}-a_{3}b_{0}c_{1}+a_{3}b_{1}c_{0},\\
R&=&a_{0}b_{1}c_{2}-a_{1}b_{0}c_{2}+a_{2}b_{0}c_{1}-a_{2}b_{1}c_{0},\\
S&=&a_{1}b_{3}c_{2}-a_{2}b_{3}c_{1}+a_{3}b_{1}c_{2}.
\end{eqnarray*}
Since all parameters of system (\ref{eq:model}) are positive, we have $K>0$, $A>0$, and $C>0$. Thus, $E_0$ and $E_1$ always exist. For $E_2$ and $E_3$ to exist, we require $B>0$ and $D>0$, respectively, or equivalently $A<K$ and $C<K$, respectively. Now, if $S$ is positive (resp. negative), then each of $P$, $Q$, and $R$ must also be positive (resp. negative) for $E_4$ to be a positive equilibrium. The following theorem summarizes these results.

\begin{Theorem} For system (\ref{eq:model}), the equilibrium solutions $E_0=(0,0,0)$ and $E_1=(K,0,0)$ always exist, while 
$E_2$ and $E_3$ exist provided $a_0b_1>a_1b_0$ and $a_0c_1>a_1c_0$, respectively, or equivalently, $A<K$ and $C<K$, respectively. The positive equilibrium solution $E_4$ exists if $P/S$, $Q/S$, and $R/S$ are all positive.
\label{thm:existence}
\end{Theorem}

\subsection{Local Stability of Equilibria}
Let $X(t)=[x(t),y(t),z(t)]^T$. Then, the linear system corresponding to (\ref{eq:model}) around an equilibrium solution $E_{*}=(x_{*},y_{*},z_{*})$ is  
$
\dot{X}(t)=M_0 X(t)+M_1 X(t-\tau)
$ 
where
$$
M_0=
\left[\begin{tabular}{ccc}
$a_{0} - a_{1}x_{*} - a_{2}y_{*} - a_{3}z_{*}$ &
$-a_{2}x_{*}$ & 
$-a_{3}x_{*}$ \\
$b_{1}y_{*}$ & 
$- b_{0} + b_{1}x_{*} - b_{3}z_{*}$ & 
$-b_{3}y_{*}$ \\
$c_{1}z_{*}$ & 
$c_{2}z_{*}$ & 
$- c_{0} + c_{1}x_{*} + c_{2}y_{*}$ 
\end{tabular}\right]
$$
and
$$
M_1=
\left[\begin{tabular}{ccc}
$-a_{1} x_{*}$ & 
$0$ & 
$0$ \\
$0$ & 
$0$ & 
$0$ \\
$0$ & 
$0$ & 
$0$ 
\end{tabular}\right],
$$
and with corresponding characteristic equation 
\begin{equation}
\det(\lambda I-M_0-M_1 e^{-\lambda\tau})=0.
\label{eq:chareqn}
\end{equation}
If all roots of (\ref{eq:chareqn}) have negative real part, then the equilibrium solution $E_*$ is locally asymptotically stable \cite{halsmith}. We can think of the roots of (\ref{eq:chareqn}) as continuous functions of the delay time $\tau$, that is, $\lambda=\lambda(\tau)$. 
When $\tau=0$, (\ref{eq:chareqn}) is a polynomial equation. In this case, the well-known Routh-Hurwitz criterion can then be utilize to provide stability conditions for $E_*$. 
As $\tau$ is increased from zero, Corollary 2.4 of \cite{ruanwei} tells us that the sum of the orders of the roots of (\ref{eq:chareqn}) in the open right half-plane can change only if a zero appears on or crosses the imaginary axis. That is, stability switch occurs at a critical value $\tau=\tau_0$ where $\lambda(\tau_0)$ is either zero or purely imaginary. The transversality condition $d(Re\ \lambda)/d\tau|_{\tau=\tau_0} >0$ implies that the eigenvalues cross from left to right. Hence, if $E_*$ is stable at $\tau=0$, then, as $\tau$ is increased, it loses its stabilty at $\tau=\tau_0$. Thus, $E_*$ is locally asymptotically stable when $\tau\in(0,\tau_0)$. 
If there are no roots of (\ref{eq:chareqn}) that cross the imaginary axis, then there are no stability switches and in this case, we have absolute stability \cite{absstab}. That is, 
$E_*$ remains stable for all delay time $\tau>0$.

\section{Main Results}
In the following, we give conditions for the stability of each equilibrium solution of system (\ref{eq:model}) when $\tau>0$ using the technique mentioned above. When $\tau=0$, the stability analysis of each of the equilibrium solution of system (\ref{eq:model}) can be found in \cite[Appendix A]{namba2008} and in \cite{holtpolis} but in slightly different form. We mention them here for completeness.\\


At $E_0$, the characteristic equation (\ref{eq:chareqn}) becomes 
$
(\lambda-a_0)(\lambda+b_0)(\lambda+c_0)=0
$ 
whose roots are $a_0>0$, $-b_0<$ and $-c_0<0$.
Thus, we have the following result. 

\begin{Theorem}
The equilibrium solution $E_0=(0,0,0)$ of system (\ref{eq:model}) is a local saddle point and is unstable for all delay time $\tau>0$.
\end{Theorem}

At $E_1$, the characteristic equation (\ref{eq:chareqn}) becomes
\begin{equation}
\left(\lambda-b_1 a_2 B/a_1\right)
\left(\lambda-c_1 a_3 D/a_1\right)
\left(\lambda+a_{0} e^{-\lambda\tau}\right)
=0.
\label{eq:chareqn1}
\end{equation} 
When $\tau=0$, (\ref{eq:chareqn1}) has roots $-a_0<0$, $b_1 a_2 B/a_1$, and $c_1 a_3 D/a_1$. In this case, we see that $E_1$ is locally asymptotically stable provided both $B$ and $D$ are negative, or equivalently $A>K$ and $C>K$. Recall from Theorem \ref{thm:existence} that these conditions imply that both $E_2$ and $E_3$ do not exist. Suppose now that $\tau>0$, and $A>K$ and $C>K$. It is known that all roots of $\lambda+a_{0} e^{-\lambda\tau}=0$ have negative real part if and only if $a_0\tau<\pi/2$ (see for example \cite[pp.70]{yangkuang}). 
Thus, $E_1$ is locally asymptotically stable if and only if $\tau<\pi/2a_0$. 

\begin{Theorem}
\label{thm:e1}
Suppose that (\ref{eq:model}) satisfies $a_1b_0>a_0b_1$ and $a_1c_0>a_0c_1$, or equivalently $A>K$ and $C>K$, respectively. Then, the equilibrium solution $E_1=(K,0,0)$ is locally asymptotically stable if and only if $\tau\in(0,\tau_c)$ where $\tau_c=\pi/2a_0$. If $\tau=\tau_c$, then (\ref{eq:model}) undergoes a Hopf bifurcation at $E_1$.
\end{Theorem}

\begin{figure}[H]
	\centering
		\includegraphics[width=8.00cm]{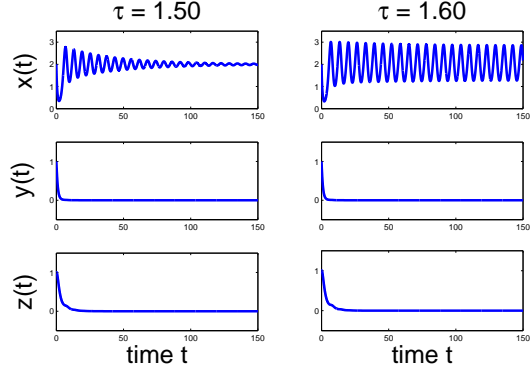}
	\caption{Stability switch occurs at a Hopf bifurcation point where $\tau=\pi/2$. $E_1$ is stable when $\tau<\pi/2$ (Left) and is unstable when $\tau>\pi/2$ (Right).}
\label{fig:e1}
\end{figure}

\noindent\textbf{Example 1.} The set of parameter values
$$
(a_0,a_1,a_2,a_3,b_0,b_1,b_3,c_0,c_1,c_2)=(1,0.5,1,0.6,0.75,0.25,0.5,0.5,0.15,0.3)
$$ 
satisfies the assumptions in Theorem \ref{thm:e1}, and gives $\tau_c=\pi/2$ (approx. 1.5708). Figure \ref{fig:e1} illustrates the stability switch at $\tau=\pi/2$ using the initial condition $(x(t),y(t),z(t))=(2,1,1)$ for $t\in[-\tau,0]$. $E_1=(2,0,0)$ is stable for $\tau\in(0,\tau_c)$.

At $E_2$, the characteristic equation (\ref{eq:chareqn}) becomes
\begin{equation}
\left(\lambda-R/a_2 b_1\right)
\left(\lambda^{2}+\bar{b}\lambda e^{-\lambda\tau}+\bar{c}\right)
=0
\label{eq:chareqne2}
\end{equation}
where $\bar{b}=a_1 A$ and $\bar{c}=a_2 b_0 B$.
Since the existence of $E_2$ requires both $A$ and $B$ to be positive, we have $\bar{b}>0$ and $\bar{c}>0$. 
First, consider the case when $\tau=0$.  
Since $\bar{b}>0$ and $\bar{c}>0$, the Routh-Hurwitz criterion tells us that all roots of 
$
\lambda^{2}+\bar{b}\lambda +\bar{c}=0
$ 
have negative real parts. Thus, $E_2$ is locally asymptotically stable if and only if $R<0$. 
Suppose now that $\tau>0$ and $R<0$. 
Since both $\bar{c}$ and $R$ are non-zero, then $\lambda=0$ is not a root of (\ref{eq:chareqne2}). Meanwhile, notice that $\lambda=i\mu$, with $\mu>0$, is a root of (\ref{eq:chareqne2}) if and only if $\mu$ satisfies  
$
-\mu^2+i\bar{b}\mu e^{-i\mu\tau}+\bar{c}=0.
$ 
Splitting into real and imaginary parts, we obtain 
\begin{equation}
\bar{b}\mu\sin(\mu\tau) = \mu^2 - \bar{c}
\qquad\mbox{and}\qquad
\bar{b}\mu\cos(\mu\tau) = 0.
\label{eq:sincose2}
\end{equation}
Squaring each sides of equations in (\ref{eq:sincose2}) and then adding them together, we get
$
\mu^4-\left(\bar{b}^2 + 2\bar{c}\right)\mu^2 +\bar{c}^2 =0,
$  
or equivalently, 
$
\mu^2
=
\frac{1}{2}\left[
\bar{b}^2 + 2\bar{c}
\pm\sqrt{\left(\bar{b}^2 +2\bar{c}\right)^2-4\bar{c}^2}\right]
$. 
Since $\bar{b}^2 + 2\bar{c} >0$ and the discriminant 
$
(\bar{b}^2 +2\bar{c})^2-4\bar{c}^2
=\bar{b}^2(\bar{b}^2+4\bar{c})
>0
$, we obtain two positive values of $\mu$ given by
\begin{equation}
\mu_{\pm}
=
\left\{\frac{1}{2}\left[\left(\bar{b}^2 + 2\bar{c}\right)\pm\sqrt{\left(\bar{b}^2 +2\bar{c}\right)^2-4\bar{c}^2}\right]\right\}^{1/2}
\label{eq:mupm}
\end{equation}
with $0<\mu_{-}<\mu_{+}$. Using (\ref{eq:sincose2}), we get $\cot\mu\tau=0$. So that $\mu\tau=(2k+1)\pi/2$ for $k=0,1,2,...$. Now, corresponding to a given $\mu$, define the increasing sequence 
$\tau_{k}=(2k+1)\pi/2\mu$ for $k=0,1,2,\dots$. The proof of the following lemma can be found in \cite[pp. 74--75]{yangkuang}.

\begin{Lemma}
Let $\lambda(\tau)$ be the root of (\ref{eq:chareqne2}) satisfying $\lambda(\tau_k)=i\mu$. Then, 
$$
\mbox{sign}\ \left.\frac{d(Re\ \lambda)}{d\tau} \right|_{\tau=\tau_k}
=
\mbox{sign}\ (2\mu^2 - \bar{b}^2-2\bar{c}).
$$
\label{lem:sign2}
\end{Lemma}
For $\mu=\mu_{\pm}$, define the sequence $\tau^{\pm}_{k}=(2k+1)\pi/2\mu_{\pm}$ for $k=0,1,2,\dots$ correspondingly, and let $\tau^{\pm}_{0}=\tau_{\mu_{\pm}}$. 
Since $0<\mu_{-}<\mu_{+}$, we have $\tau_{\mu_{+}}<\tau_{\mu_{-}}$, and hence $\tau_{\mu_+}$ is the smallest amongst all $\tau_{k}^{\pm}$. 
From (\ref{eq:mupm}), we have  
$
2\mu_{+}^2 - \bar{b}^2-2\bar{c}
>0
$. 
Thus, at $\tau=\tau_{\mu_+}$ (or equivalently at $\mu=\mu_+$), the quantity 
$
\frac{d(Re\ \lambda)}{d\tau}
$ 
is positive by Lemma \ref{lem:sign2}. 
This, together with the Hopf bifurcation theorem \cite{halsmith}, we have the following result.


\begin{Theorem}
\label{thm:e2}
Suppose that in system (\ref{eq:model}), $a_{0}b_{1}c_{2}+a_{2}b_{0}c_{1} <a_{1}b_{0}c_{2}+a_{2}b_{1}c_{0}$ or equivalently $R<0$. Then, the equilibrium solution $E_2=(A,B,0)$ of (\ref{eq:model}) is locally asymptotically stable whenever $\tau\in(0,\tau_{\mu_{+}})$. If $\tau=\tau_{\mu_+}$, then system (\ref{eq:model}) undergoes a Hopf bifurcation at $E_2$.
\end{Theorem}

\begin{figure}[H]
	\centering
		\includegraphics[width=8.00cm]{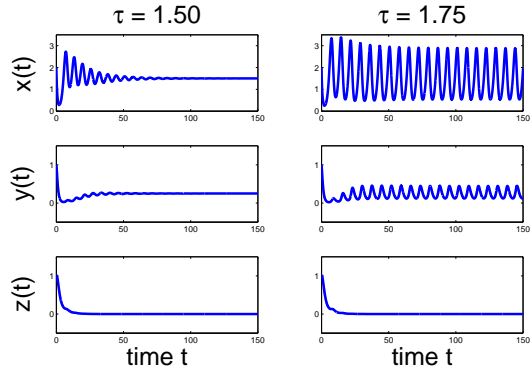}
	\caption{Stability switch occurs at a Hopf bifurcation where $\tau=\tau_{\mu_+}=1.6573$ approximately. The equlibrium solution $E_2=(1.50,0.25,0.00)$ is stable when $\tau<\tau_{\mu_+}$ (Left) and is unstable when $\tau>\tau_{\mu_+}$ (Right).}
\label{fig:e2}
\end{figure}

\noindent\textbf{Example 2.} Using the same set of parameters and initial condition as in Example 1 with $b_1$ changed to $0.50$, the assumptions in Theorem \ref{thm:e2} are satisfied. Figure \ref{fig:e2} illustrates the stability switch at $\tau=\tau_{\mu_+}=1.6573$ approximately.

A similar condition for the stability of $E_3$ can be obtained using the same analysis used for $E_2$, and is given in the following theorem.

\begin{Theorem}
Suppose that in system (\ref{eq:model}), 
$a_{1}b_{3}c_{0}+a_{3}b_{1}c_{0}
<a_{0}b_{3}c_{1}+a_{3}b_{0}c_{1}$ or equivalently $Q<0$. Let $\tilde{b}=a_1 C$, $\tilde{c}=a_3 c_0 D$, and 
$$
\nu_{+}
=
\left\{\frac{1}{2}\left[\left(\tilde{b}^2 + 2\tilde{c}\right)
+
\sqrt{\left(\tilde{b}^2 +2\tilde{c}\right)^2-4\tilde{c}^2}\right]\right\}^{1/2}.
$$
Then, the equilibrium solution $E_3=(C,0,D)$ of (\ref{eq:model}) is locally asymptotically stable whenever $\tau\in(0,\tau_{\nu_{+}})$ where $\tau_{\nu_{+}}=\pi/2\nu_{+}$. If $\tau=\tau_{\nu_+}$, then system (\ref{eq:model}) undergoes a Hopf bifurcation at $E_3$.
\end{Theorem}

\subsection*{Stability of the Positive Equilibrium $E_4$}

At $E_4$, the characteristic equation (\ref{eq:chareqn}) becomes
\begin{eqnarray}
\lambda^3+a\lambda+b+(c\lambda^2+d)e^{-\lambda\tau}=0
\label{eq:chareqnposeq}
\end{eqnarray}
where 
$a=(PQa_{2}b_{1}+PRa_{3}c_{1}+QRb_{3}c_{2})/S^2$, 
$b=PQR(a_{3}b_{1}c_{2}-a_{2}b_{3}c_{1})/S^3$,  
$c=a_{1}P/S$, and 
$d=PQRa_{1}b_{3}c_{2}/S^3.$ 
Note here that 
$
b+d
=PQR/S^2.
$ Order-three quasi-polynomials with single delay, different to (\ref{eq:chareqnposeq}), have been studied in \cite{culshawruan, katriruan}. We follow a similar method used in \cite{culshawruan} to analyze the distribution of the roots of (\ref{eq:chareqnposeq}) on the complex plane. 
When $\tau=0$, then (\ref{eq:chareqnposeq}) reduces to 
\begin{equation}
\lambda^{3}
+c\lambda^{2}
+a\lambda
+(b+d)
=0.
\label{eq:poly}
\end{equation}
The Routh-Hurwitz criterion tells us that all roots of (\ref{eq:poly}) have negative real part if and only if the following inequalities hold: 
$
a,c,b+d>0
$ 
and 
$
ac-(b+d)>0.
$ 
If $S<0$, then $P$, $Q$, and $R$ must all be negative for $E_4$ to be a positive equilibrium. This implies that $a>0$, $c<0$, and $b+d<0$, and thus, $E_4$ is unstable. If $S>0$, then each of $P$, $Q$, and $R$ must be positive for $E_4$ to be a positive equilibrium. In this case, $a$, $c$, and $b+d$ are all positive. Now, notice that if in addition to the assumption that $S>0$, we also have $b<0$ (or equivalently $a_2b_3c_1>a_3b_1c_2$), then 
$
ac-(b+d)
=a_1 P^2 (Q a_2b_1+ R a_3c_1)/S^3-b>0.
$ 
Hence, for the case when $\tau=0$, $E_4$ is locally asymptotically stable whenever $S>0$ and $b<0$.\\

Suppose now that $\tau>0$, and assume that $S>0$ and $b<0$. 
If $\lambda=0$ is a root of (\ref{eq:chareqnposeq}), then $b+d=0$. However, $b+d\neq0$ since $b+d=PQR/S^2$ and none of $P$, $Q$, and $R$ is equal to zero. Thus,  $\lambda=0$ is not a root of (\ref{eq:chareqnposeq}). If $\lambda=i\omega$, with $\omega>0$, is a root (\ref{eq:chareqnposeq}), then 
$
-i\omega^3+ia\omega+b+(-c\omega^2+d)e^{-i\omega\tau}=0.
$ 
Splitting into real and imaginary parts, we get
\begin{equation}
(c\omega^2-d)\cos\omega\tau
=b,
\qquad\mbox{and}\qquad
(c\omega^2-d)\sin\omega\tau
=\omega^3-a\omega.
\label{eq:sincos}
\end{equation}
Squaring each sides of these equations and then adding corresponding sides gives
\begin{equation}
\omega^6+\alpha\omega^4+\beta\omega^2+\gamma=0
\label{eq:deg6}
\end{equation}
where 
$\alpha = -2a-c^2$, 
$\beta  = a^2+2cd$, and 
$\gamma =b^2-d^2$. 
We claim that $\alpha<0$, $\beta>0$, and $\gamma<0$. 
To see this, observe that since $S>0$, we have $P,Q,R>0$ and consequently, $a$, $c$, and $d$ are all positive. Immediately, we see that $\alpha<0$ and $\beta>0$. 
Meanwhile, the assumptions $S>0$ and $b<0$ gives $b+d>0$ and $b-d<0$, respectively. As a result, $\gamma=(b+d)(b-d)<0$. We now show that (\ref{eq:deg6}) has at least one positive root. Let $h(u)=u^3+\alpha u^2+\beta u+\gamma$, and note that if the equation $h(u)=0$ has a positive root $u=u_0$, then (\ref{eq:deg6}) has a positive root $\omega_0=\sqrt{u_0}$.
Now, notice that 
$
h'(u)=3u^2+2\alpha u+\beta=0
$ 
has two positive roots  
$
u_{\pm}=(-\alpha\pm\sqrt{\alpha^2-3\beta}\ )/3
$ 
with $0<u_{-}<u_{+}$. Since $h(0)=\gamma<0$, $h'(0)=\beta>0$, and $h''(0)=2\alpha<0$, we know that the graph of $h$ passes through the point $(0,\gamma)$ below the horizontal axis, and then increases in a concave down manner. The continuity of $h$ and the fact that $h$ increases without bound as $u\rightarrow+\infty$ guarantee that $h$ has at least one positive root $u_0$. That is, (\ref{eq:deg6}) has a positive root $\omega_0=\sqrt{u_0}$, and therefore (\ref{eq:chareqnposeq}) has purely imaginary roots $\lambda=\pm i\omega_0$. Noting that a local maximun of $h$ occurs at $u=u_{-}$ and a local minimun of $h$ occurs at $u=u_{+}$, we see that the graph of the cubic polynomial $h$ is increasing on $(0,u_{-})$ and on $(u_{+},+\infty)$. Moreover, we either have $0<u_0<u_{-}$ or $u_0>u_{+}$.\\

If the equation $h(u)=0$ has more than one positive roots, then it must have exactly three positive roots so that $\omega$ is a simple root of (\ref{eq:deg6}). 
Using the first equation in (\ref{eq:sincos}), for a given  $\omega$, define its corresponding increasing sequence 
$\tau_k=\frac{1}{\omega}\left[\cos^{-1}\left(\frac{b}{c\omega^2-d}\right)+2\pi k\right]$ for $k=0,1,2,...$. Specifically, for a given $\omega$, we get a corresponding $\tau_0$. 
Among the three positive roots of the equation $h(u)=0$, we then choose $u_0$ so that $\omega=\omega_0=\sqrt{u_0}$ has corresponding $\tau_0$ that is smallest. This guarantees that a pair of purely imaginary eigenvalues will first occur at $\tau=\tau_0$.

\begin{Lemma} 
Let $\lambda(\tau)$ be the root of (\ref{eq:chareqnposeq}) satisfying $\lambda(\tau_k)=i\omega_0$. Then, 
$$
\left.\frac{d(Re\ \lambda)}{d\tau}\right|_{\tau=\tau_k}>0.
$$
\label{lem:signposeq}
\end{Lemma}
\textbf{Proof.} From (\ref{eq:chareqnposeq}), 
$$
\left[(3\lambda^2+a)+(2c\lambda-(c\lambda^2+d)\tau)e^{-\lambda\tau}\right]\frac{d\lambda}{d\tau}
-\lambda(c\lambda^2+d)e^{-\lambda\tau}=0.
$$
Consequently,
\begin{eqnarray*}
\left(\frac{d\lambda}{d\tau}\right)^{-1}
&=&
\frac{2c\lambda+(3\lambda^2+a)e^{\lambda\tau}}{\lambda(c\lambda^2+d)}
-\frac{\tau}{\lambda}\\
&=&
\frac{2c}{c\lambda^2+d}
-\frac{3\lambda^2+a}{\lambda(\lambda^3+a\lambda+b)}
-\frac{\tau}{\lambda}
\end{eqnarray*}
since 
$e^{\lambda\tau}/(c\lambda^2+d)
=
-1/(\lambda^3+a\lambda+b)
$ using (\ref{eq:chareqnposeq}). Hence,
\begin{eqnarray*}
\mbox{sign} \left\{\frac{d\left(Re\ \lambda\right)}{d\tau}\right\}_{\lambda=i\omega_0}
&=&
\mbox{sign} \left\{Re \left(\frac{d\lambda}{d\tau}\right)^{-1}\right\}_{\lambda=i\omega_0}\\
&=&
\mbox{sign} \left\{Re\ \frac{2c}{c\lambda^2+d}
+Re\ \frac{-3\lambda^2-a}{\lambda^4+a\lambda^2+b\lambda}
\right\}_{\lambda=i\omega_0}\\
&=&
\mbox{sign} \left\{Re\ \frac{2c}{-c\omega_0^2+d}
+Re\ \frac{3\omega_0^2-a}{\omega_0^4-a\omega_0^2+ib\omega_0}
\right\}\\
&=&
\mbox{sign} \left\{\frac{-2c}{c\omega_0^2-d}
+\frac{(3\omega_0^2-a)(\omega_0^4-a\omega_0^2)}{(\omega_0^4-a\omega_0^2)^2+(b\omega_0)^2}
\right\}\\
&=&
\mbox{sign} \left\{\frac{-2c}{c\omega_0^2-d}
+\frac{(3\omega_0^2-a)(\omega_0^2-a)}{(\omega_0^3-a\omega_0)^2+b^2}
\right\}\\
&=&
\mbox{sign} \left\{\frac{-2c}{c\omega_0^2-d}
+\frac{(3\omega_0^2-a)(\omega_0^2-a)}{(c\omega_0^2-d)^2}
\right\}
\end{eqnarray*}
since $(\omega_0^3-a\omega_0)^2+b^2=(c\omega_0^2-d)^2$ from (\ref{eq:sincos}). Thus,
\begin{eqnarray*}
\mbox{sign} \left\{
\frac{d\left(Re\ \lambda\right)}{d\tau}\right\}_{\lambda=i\omega_0}
&=&
\mbox{sign} \left\{
-2c(c\omega_0^2-d)+(3\omega_0^2-a)(\omega_0^2-a)
\right\}\\
&=&
\mbox{sign} \left\{
3\omega_0^4+(-4a-2c^2)\omega_0^2+(a^2+2cd)
\right\}\\
&=&
\mbox{sign} \left\{3\omega_0^4+2\alpha\omega_0^2+\beta\right\}.
\end{eqnarray*}
Since $h$ is increasing on the intervals $(0,u_{-})$ and $(u_{+},\infty)$, and $u_0$ belongs to $(0,u_{-})$ or to $(u_{+},\infty)$, we have $h'(u_0)=3u_0^2+2\alpha u_0+\beta>0$. 
This implies that $3\omega_0^4+2\alpha\omega_0^2+\beta>0$ since $\omega_0=\sqrt{u_0}$. Thus,
$$
\mbox{sign} \left\{
\frac{d\left(Re\ \lambda\right)}{d\tau}\right\}_{\lambda=i\omega_0}
=
\mbox{sign} \left\{3\omega_0^4+2\alpha\omega_0^2+\beta\right\}
=+1
$$
and this completes the proof. \\

By Lemma \ref{lem:signposeq}, 
$
\left.\frac{d(Re\ \lambda)}{d\tau}\right|_{\tau=\tau_0}>0.
$ 
This, together with the Hopf bifurcation theorem \cite{halsmith} give the following result.

\begin{Theorem}
\label{thm:e4}
Suppose that in system (\ref{eq:model}), 
$
a_{1}b_{3}c_{2}+a_{3}b_{1}c_{2}>a_{2}b_{3}c_{1} 
$ 
and 
$
a_2b_3c_1>a_3b_1c_2
$ 
(or equivalently $S>0$ and $b<0$, respectively). 
Then, the positive equilibrium solution $E_4$ of (\ref{eq:model}) is locally asymptotically stable whenever $\tau\in(0,\tau_0)$. If $\tau=\tau_0$, then system (1) undergoes a Hopf bifurcation at $E_4$.
\end{Theorem}

\begin{figure}[H]
	\centering
		\includegraphics[height=3.75cm]{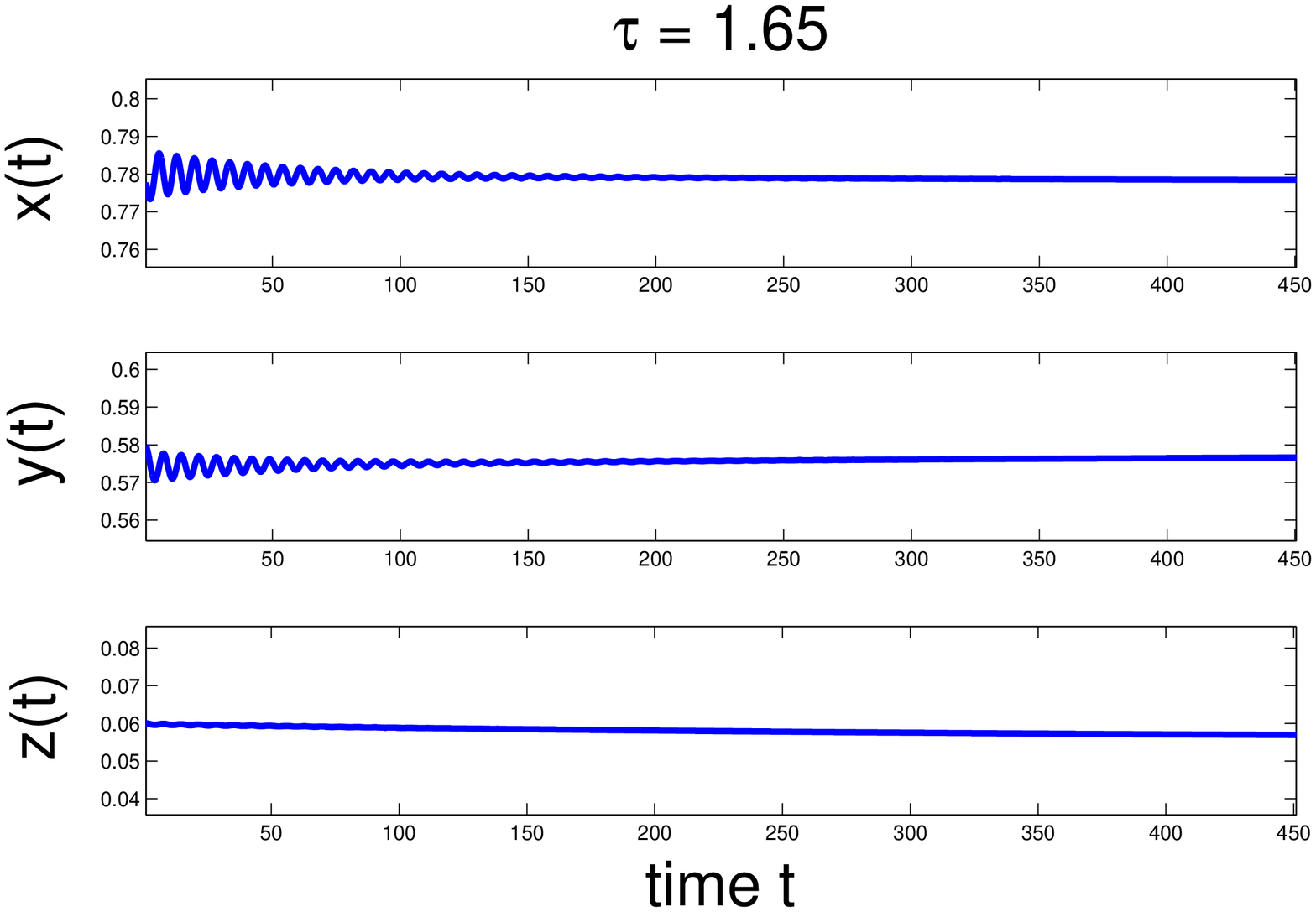}\includegraphics[height=3.75cm]{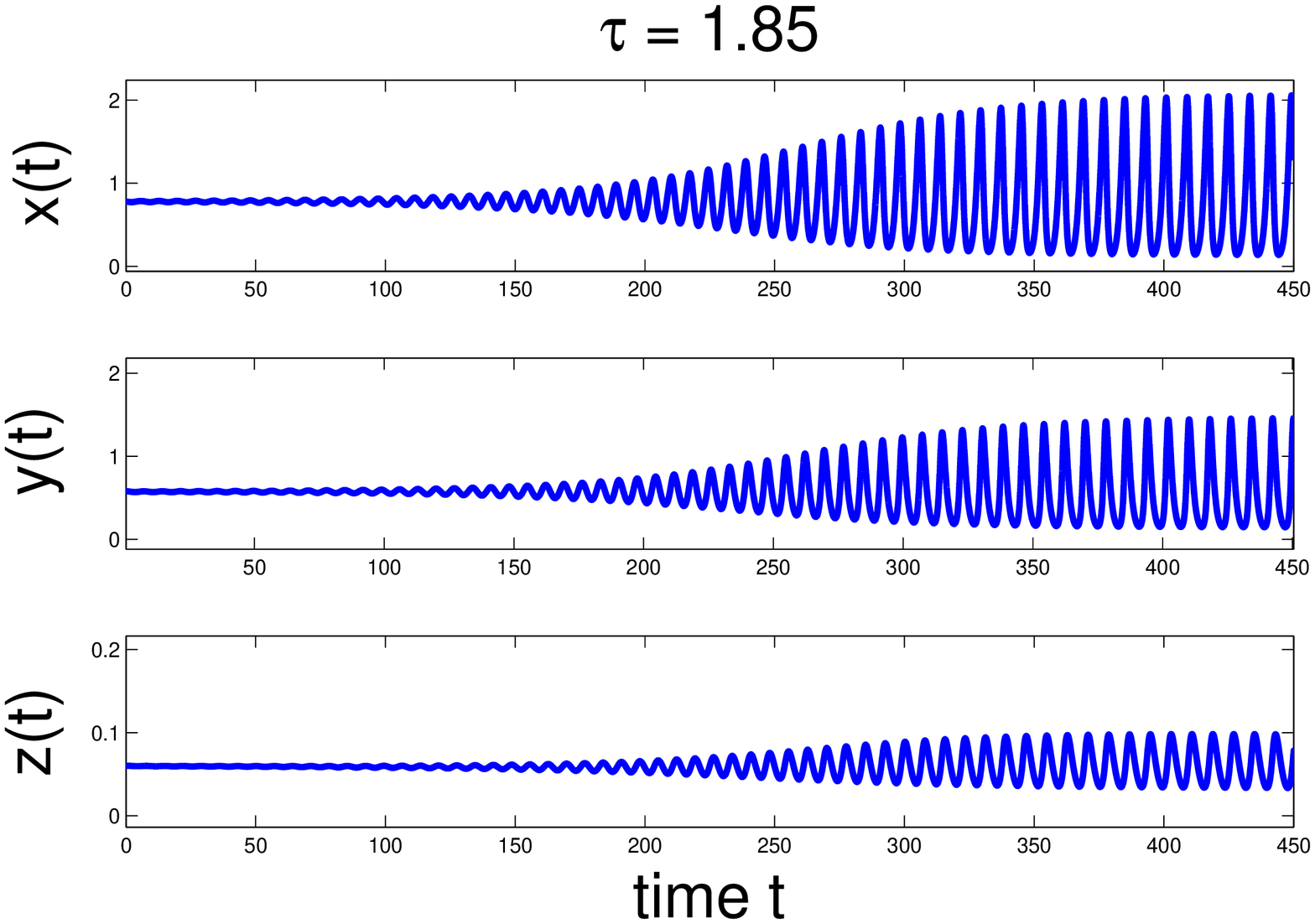}
	\caption{Stability switch occurs at $\tau=\tau_0=1.7438$ (approx.). The positive equilibrium solution $E_4=(0.7778,0.5778,0.0556)$ is stable for $\tau<\tau_0$ (Left) and is unstable when $\tau>\tau_0$ (Right).}
\label{fig:e4}
\end{figure}

\noindent\textbf{Example 3.} Using the same set of parameters as in Example 1 with $b_1=1$ and $c_1=0.42$, the assumptions in Theorem \ref{thm:e4} are satisfied. Using the initial condition $(x(t),y(t),z(t))=(0.78, 0.58, 0.06)$ for $t\in[-\tau,0]$, Figure \ref{fig:e4} illustrates the stability switch at $\tau=\tau_0=1.7438$ approximately. 

\section{Numerical Continuation}
\begin{figure}[H]
	\centering
		\includegraphics[width=7.00cm]{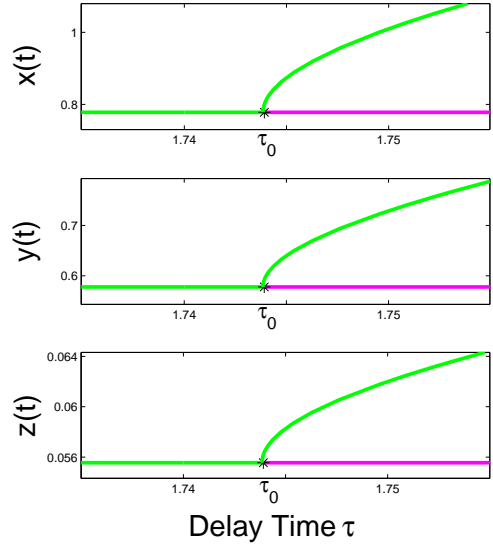}
	\caption{Stable branch of periodic solutions emerges from the Hopf bifurcation point marked with ($\ast$) where $\tau=\tau_0$.}
\label{fig:e4bifdiag}
\end{figure}

Recall that for $\tau=0$, $E_4$ is locally asymptotically stable whenever $S>0$ and $b<0$. 
A branch of equilibrium solutions can be obtained by following or continuing this equilibrium solution in DDE-Biftool \cite{ddebiftoolarticle} using the delay time $\tau$ as parameter. DDE-Biftool is a numerical continuation and bifurcation analysis tool developed by Engelborghs et al \cite{ddebiftoolmanual}. We use the same set of parameters and initial condition used in Example 3. In Figure \ref{fig:e4bifdiag}, the top panel shows the branch of equilibrium solutions (horizontal line) obtained in DDE-Biftool, where green and magenta represents the stable and unstable parts of the branch, respectively. A change of stability occurs at a Hopf bifurcation point marked with ($\ast$) where $\tau=\tau_0$. Again, we use DDE-biftool to continue this Hopf point into a branch of periodic solutions. We obtained a stable branch of periodic solutions and this shown as the green curve in the top panel of Figure \ref{fig:e4bifdiag}. Here, the vertical axis gives a measure of the maximum value of the oscillation of the periodic solutions in the $x(t)$ component. The middle and bottom panels corresponds to the same bifurcation diagram as the top panel but for the $y(t)$ and the $z(t)$ components.\\

Figure \ref{fig:e4bifdiag} illustrates the results of Theorem \ref{thm:e4}, that is, the stability switch at $\tau=\tau_0$ and the occurrence of Hopf bifurcation. Numerical continuations showed that, as the positive equilibrium $E_4$ loses its stability, the branch of periodic solutions that emerges from the Hopf point is stable. Biologically, this means that the three species still persist and the only difference now is that the three populations are oscillating.   

\section{Conclusion}
We studied a delayed Lotka-Volterra IG predation model where, in the absence of predation, the basal resource grows according to the delayed logistic equation. By first considering the case where $\tau=0$, stability conditions for the all non-negative equilibria are established using the well-known Routh-Hurwitz criterion. For $\tau>0$, we showed that stability switch occurs at a Hopf bifurcation where the stable equilibrium becomes unstable. Moreover, numerical continuation shows that a stable branch of periodic solutions emerges from the Hopf point as the positive equilibrium solution becomes unstable. This means that by increasing the delay time the positive equilibrium could become unstable. However, since the periodic orbit that will emerge is stable, all three species will still persist.

\end{document}